\documentclass{amsart}

\usepackage{color}
\usepackage[dvipsnames]{xcolor}

\usepackage{xargs}
\usepackage{expl3}

\usepackage[utf8]{inputenc}				
\usepackage[english]{babel}
\usepackage[T1]{fontenc}				

\usepackage{stix2}
\usepackage{amsmath,amsthm,amsfonts}
\usepackage{mathtools}
\usepackage{bm}
\usepackage{nicefrac}
\usepackage{extarrows}
\usepackage{leftidx}

\usepackage{wrapfig}
\usepackage{float}
\usepackage{subcaption}
\usepackage{epsfig}
\usepackage{epstopdf}
\usepackage{svg}
\usepackage{pdfpages}
\usepackage{graphicx}
\usepackage{tikz}
\usetikzlibrary{math,arrows,intersections}
\usepackage{transparent} 
\usepackage{pgfplots}
\pgfplotsset{compat = newest}
\pgfplotstableset{col sep=comma}

\usepackage{listings}
\usepackage{algorithm}
\usepackage{algpseudocode}
\usepackage[framed,numbered,autolinebreaks,useliterate]{mcode}

\usepackage{ltxtable}
\usepackage{tabularx}
\usepackage{multirow}
\usepackage{longtable}
\usepackage{booktabs}

\usepackage{paralist}

\usepackage{url}
\usepackage[colorlinks = false,
			pdfpagelabels,
			pdfstartview = FitH,
			bookmarksopen = true,
			bookmarksnumbered = true,
			linkcolor = black,
			plainpages = false,
			hypertexnames = false,
			citecolor = black,
			hidelinks] {hyperref}

\usepackage[
backend=biber,
style=numeric,
sorting=nyt,
giveninits=true,
doi=false,
isbn=false,
url=false,
maxbibnames=99,
language=english,
]{biblatex}
\renewbibmacro{in:}{}
\usepackage{cleveref}
\newcommand{\crefl}[1]{\labelcref{#1}}

\usepackage{lipsum}
\usepackage{changes}
\usepackage{csquotes}

\makeatletter
\newcommand\RedeclareMathOperator{%
  \@ifstar{\def\rmo@s{m}\rmo@redeclare}{\def\rmo@s{o}\rmo@redeclare}%
}
\newcommand\rmo@redeclare[2]{%
  \begingroup \escapechar\m@ne\xdef\@gtempa{{\string#1}}\endgroup
  \expandafter\@ifundefined\@gtempa
     {\@latex@error{\noexpand#1undefined}\@ehc}%
     \relax
  \expandafter\rmo@declmathop\rmo@s{#1}{#2}}
\newcommand\rmo@declmathop[3]{%
  \DeclareRobustCommand{#2}{\qopname\newmcodes@#1{#3}}%
}
\@onlypreamble\RedeclareMathOperator
\makeatother
\newcommand{\RedeclarePairedDelimiter}[3]{
	\let\angles\relax\DeclarePairedDelimiter{#1}{#2}{#3}}
	
\newcommand	{\ts} [1]{{\textstyle{#1}}} 	

\newcommand {\ra}	{\rightarrow}

\newcommand {\R}	{\mathbb{R}}
\newcommand {\N}	{\mathbb{N}}
\newcommand	{\C}	{\mathbb{C}}

\newcommand {\cP}   {\mathcal{P}}
\newcommand {\cN}   {\mathcal{N}}

\newcommand	{\bintegral}[4]	{\int_{#1}^{#2} \! #3 \, \mathrm{d}#4}
\DeclarePairedDelimiter  {\norm} {\lVert} {\rVert}
\DeclarePairedDelimiter  {\abs}  {\lvert} {\rvert}
\DeclarePairedDelimiter  {\opnorm}{\Vvert} {\Vvert}

\DeclarePairedDelimiter  {\inner}{(}{)}

\DeclarePairedDelimiter  {\brackets}{\lbrack}{\rbrack}
\DeclarePairedDelimiter  {\braces}  {\lbrace}{\rbrace}
\DeclarePairedDelimiter  {\parens}  {\lparen}{\rparen}
\RedeclarePairedDelimiter{\angles}  {\langle}{\rangle}

\DeclareMathOperator  {\spann} {span}

\RedeclareMathOperator{\closure}{clos}
\RedeclareMathOperator{\div}   {div}
\DeclareMathOperator  {\kernel}{ker}
\DeclareMathOperator  {\ind}{ind}
\makeatletter
\renewcommand*\env@matrix[1][*\c@MaxMatrixCols c]{%
  \hskip -\arraycolsep
  \let\@ifnextchar\new@ifnextchar
  \array{#1}}
\makeatother

\makeatletter
\newsavebox\myboxA
\newsavebox\myboxB
\newlength\mylenA

\newcommand*\xoverline[2][0.75]{%
    \sbox{\myboxA}{$\m@th#2$}%
    \setbox\myboxB\null
    \ht\myboxB=0.9\ht\myboxA%
    \dp\myboxB=\dp\myboxA%
    \wd\myboxB=#1\wd\myboxA
    \sbox\myboxB{$\m@th\overline{\copy\myboxB}$}
    \setlength\mylenA{\the\wd\myboxA}
    \addtolength\mylenA{-\the\wd\myboxB}%
    \ifdim\wd\myboxB<\wd\myboxA%
       \rlap{\hskip 0.5\mylenA\usebox\myboxB}{\usebox\myboxA}%
    \else
        \hskip -0.5\mylenA\rlap{\usebox\myboxA}{\hskip 0.5\mylenA\usebox\myboxB}%
    \fi}
\makeatother
\graphicspath{{Images/}}

\definechangesauthor[name=Emil,  color=Mahogany]{EB}
\definechangesauthor[name=Moritz,color=NavyBlue]{MF}
\definechangesauthor[name=Niklas,color=YellowOrange]{NR}
\definechangesauthor[name=Urban, color=ForestGreen]{KU}
\definechangesauthor[name=Correction Moritz, color=purple]{MF2}

\newtheorem{theorem}{Theorem}[section]

\newtheorem{lemma}[theorem]{Lemma}

\theoremstyle{definition}

\theoremstyle{remark}
\newtheorem{remark}[theorem]{Remark}

\numberwithin{equation}{section}

\pgfplotsset{
  log x ticks with fixed point/.style={
      xticklabel={
        \pgfkeys{/pgf/fpu=true}
        \pgfmathparse{exp(\tick)}%
        \pgfmathprintnumber[fixed relative, precision=3]{\pgfmathresult}
        \pgfkeys{/pgf/fpu=false}
      }
  },
  log y ticks with fixed point/.style={
      yticklabel={
        \pgfkeys{/pgf/fpu=true}
        \pgfmathparse{exp(\tick)}%
        \pgfmathprintnumber[fixed relative, precision=3]{\pgfmathresult}
        \pgfkeys{/pgf/fpu=false}
      }
  }
}

\begin{document}

\author[Beurer]{Emil Beurer}
\address{Ulm University, Institute for Numerical Mathematics, Helmholtzstr.\ 20, D-89081 Ulm, Germany}
\curraddr{}
\email{\{emil.beurer,moritz.feuerle,niklas.reich,karsten.urban\}@uni-ulm.de}
\thanks{}

\author[Feuerle]{Moritz Feuerle}
\thanks{}

\author[Reich]{Niklas Reich}
\thanks{}

\author[Urban]{Karsten Urban}
\thanks{}

\title[An ultraweak variational method for linear PDAEs]{An ultraweak variational method for parameterized linear differential-algebraic equations}%
\thanks{The authors have no competing interests to declare that are relevant to the content of this article. We acknowledge support by the state of Baden-W\"urttemberg through bwHPC}

\date{}
\dedicatory{}

\begin{abstract}
    We investigate an ultraweak variational formulation for (parameterized) linear differential-algebraic equations (DAEs) w.r.t.\ the time variable which yields an optimally stable system. This is used within a Petrov-Galerkin method to derive a certified detailed discretization which provides an approximate solution in an ultraweak setting as well as for model reduction w.r.t.\ time in the spirit of the Reduced Basis Method (RBM). A computable sharp error bound is derived. Numerical experiments are presented that show that this method yields a significant reduction and can be combined with well-known system theoretic methods such as Balanced Truncation to reduce the size of the DAE.
\end{abstract}

\maketitle

\section{Introduction}\label{Sec:1}

Differential-Algebraic Equations (DAEs) are widely used to model several processes in science, engineering, medicine and other fields. Theory and numerical approximation methods have intensively been studied in the literature, see e.g.\ 
\cite{%
brenan1995numerical,
HairerWanner2,
KunkelMehrmann,%
Trenn2013
}, or
\cite{DAE2,DAE1}, 
which are the first two books in a forum series on DAEs. Quite often, the dimension of DAEs modeling realistic problems is so large that an efficient numerical solution (in particular in realtime environments or within optimal control) is impossible. To address this issue, Model Order Reduction (MOR) techniques have been developed and successfully applied. There is a huge amount of literature, we just mention 
\cite{
BennerSIAM2017,
Benner2017,
MR4171473,%
MehrmannStykel2005,%
Stykel2003,
Stykel2004
}. 

All methods described in those references address a reduction of the dimension of the system, whereas the temporal discretization is untouched. This paper starts at this point. We have been working on space-time variational formulations for (parameterized) partial differential equations (PPDEs) over the last decade. One particular issue has been the stability of the arising discretization which admits tight error-residual relations and thus builds the backbone for model reduction. It turns out that an \emph{ultraweak} formulation is the right tool to achieve this goal. In \cite{Feuerle2019}, we have used this framework for deriving an optimally stable variational formulation of linear time-invariant systems (LTIs). In this paper, we extend the ultraweak framework to (parameterized) DAEs and show that this can be combined with system theoretic methods such as Balanced Truncation (BT, \cite{Stykel2003}) to derive a reduction in the system dimension and time discretization size.

\subsection{Differential-algebraic equations (DAEs)}\label{Sec:DAE}
Let $E,A\in \R^{n\times n},\ n\in\N$, be two matrices ($E$ is typically singular), $I = (0, T)$, $T>0$, a time interval, $x_0\in\R^n$ some initial value and $f: I\ra\R^n$ a given right-hand side. Then, we are interested in the solution $x:I\ra\R^n$ (the state) of the following initial value problem of a linear differential-algebraic equation (DAE) with constant coefficients
\begin{align*}
        E\dot x(t) - Ax(t) &= f(t),\quad\forall t\in I,
        \qquad 
        x(0) = x_0.
\end{align*}
In order to ensure well-posedness (in an appropriate manner), we shall always assume that the initial value  $x_0$ is \emph{consistent} with the right-hand side $f$, which means that there exists some $\hat x_0\in\R^n$ such that $E\hat x_0 - Ax_0 = \lim_{t\ra 0+}f(t)$ holds. Finally, we assume that the matrix pencil 
$\braces{E,-A}$ is regular (i.e. $\det(\lambda E - A) \neq 0$ for some $\lambda\in\R$) with index $\ind\braces{E,-A} =: k \in\N$, \cite{GearPetzold83}.\footnote{Each regular matrix pencil can be transformed into \emph{Weierstrass-Kronecker canonical form}
$P(\lambda E-A)Q= \textrm{diag}(\lambda I\!d - W, \lambda N-I\!d)$ with regular matrices $P,Q\in\mathbb{C}^{n\times n}$,  \cite{Gantmacher}. The index of a regular matrix pencil $\braces{E,-A}$ is then defined by $\ind\braces{E,-A}:=\ind\braces{N}:=\min\{k\in\N: N^k=0\}$.
}

\subsection{Parameterized DAEs (PDAEs)}\label{subset:DAE}
We are particularly interested in the situation, where one does not only have to solve the above DAE once, but several times and highly efficient (e.g.\ in realtime, optimal control or cold computing devices) for different data. In order to describe that situation, we are considering a \emph{parameterized} DAE (PDAE) as follows. For some parameter vector $\mu\in\cP$, $\cP\subset\R^P$ being a compact set, we are seeking $x_\mu:I\ra\R^n$ such that 
\begin{align}\label{eq:PDAE}
        E\,\dot x_\mu(t) - A_\mu\,x_\mu(t) &= f_\mu(t),\quad\forall t\in I,
        \qquad 
        x_\mu(0) = x_{0,\mu},
\end{align}
where $A_\mu$, $f_\mu$ and $x_{0,\mu}$ are a parameter-dependent matrix, a right-hand side and an initial condition, respectively, whereas $E$ is assumed to be independent of $\mu$, see below. 
In order to be able to solve such a PDAE highly efficient for many parameters, it is quite standard to assume that parameters and variables can be separated, see e.g.\ \cite{rozza2016}. This is done by assuming a so-called \emph{affine decomposition} of the data, i.e., $E$ is (for simplicity of exposition) assumed to be parameter-independent and
\begin{align}\label{eq:affinedecomp}
    A_\mu = \sum_{q=1}^{Q_A} \vartheta^A_q(\mu)\, \tilde A_q,
    &&
    f_\mu(t) = \sum_{q=1}^{Q_f} \vartheta^f_q(\mu)\, \tilde f_q(t),
    &&
    x_{0,\mu} = \sum_{q=1}^{Q_x} \vartheta^x_q(\mu)\, \tilde x_{0,q}.
\end{align}
If such a decomposition is not given, we may produce an affinely decomposed approximation by means of the \emph{(Discrete) Empirical Interpolation Method} (\emph{(D)EIM}, \cite{Barrault2004667,DEIM}; see also \cite{MR4171473} for a system theoretic MOR for such PDAEs). For well-posedness, we assume that the matrix pencil $\braces{E,-A_\mu}$ is regular with index $\ind\braces{E,-A_\mu} = k_\mu$ for all $\mu\in\mathcal{P}$.

\subsection{Reduction to homogeneous initial conditions}\label{Sec:Hom}
Using some standard arguments, \crefl{eq:PDAE} can be reduced to homogeneous initial conditions $x_\mu(0) = 0$. To this end, construct some smooth extension of the initial data $\bar{x}_\mu\in C^1(\bar{I})^n$,  $\bar{x}_\mu(0) = x_{0,\mu}$. Then, let $\hat{x}_\mu: I \ra \R^n$ solve \crefl{eq:PDAE} with $f_\mu$ replaced by $\hat{f}_\mu := f_\mu - E\dot{\bar{x}}_\mu + A_\mu\bar{x}_\mu$ and homogeneous initial condition $\hat{x}_\mu(0) = 0$. Then, $x_\mu := \hat{x}_\mu + \bar{x}_\mu$ solves the original problem \crefl{eq:PDAE}. 
If the PDAE possess an affine decomposition \eqref{eq:affinedecomp}, it is readily seen that the modified right-hand side $\hat{f}_\mu$ also admits an affine decomposition. Hence, we can always restrict ourselves to the case of homogeneous initial conditions $x_\mu(0) = 0$, keeping in mind that variable initial conditions can be realized by different right-hand sides. 

\subsection{Organization of the material}
The remainder of this paper is organized as follows. In Section \ref{Sec:2}, we derive an ultraweak variational formulation of \eqref{eq:PDAE} and prove its well-posedness. Section \ref{Sec:3} is devoted to a corresponding Petrov-Galerkin discretization and the numerical solution, which is then used in Section \ref{Sec:MOR} to derive a certified reduced model. In Section \ref{Sec:4}, we report results of our numerical experiments and end with conclusions and an outlook in Section \ref{Sec:6}.

\section{An ultraweak variational formulation}\label{Sec:2}
It is well-known that, for any fixed parameter $\mu\in\cP$, the problem \eqref{eq:PDAE} admits a unique classical solution $x_\mu\in C^{k_\mu}(\bar I)^n$ for consistent initial conditions provided that $f_\mu\in C^{k_\mu-1}(\bar I)^n$, e.g. \cite[Lemma 2.8.]{KunkelMehrmann}. This is a severe regularity assumption, which is one of the reasons why we are interested in a variational formulation. As we shall see, an \emph{ultraweak} setting is appropriate in order to prove well-posedness, in particular stability. It turns out that this setting is also particularly useful for model reduction of \eqref{eq:PDAE} w.r.t.\ the time variable in the spirit of the reduced basis method, see \S\crefl{Sec:MOR} below.

\subsection{Ultraweak formulation of PPDEs} \label{ch:weak} 
In order to describe an ultraweak variational formulation for the above PDAE, we will review such formulations for parametric \emph{partial} differential equations (PPDEs). In particular, we are going to follow \cite{Dahmen2011} in which well-posed (ultraweak) variational forms for transport problems have been introduced, see also \cite{BSU19,Henning2021,UrbanCIME2022}. We will then transfer this framework to PDAEs in \S\crefl{Sec:vwPDAE}.

Let $\Omega\subset\R^n$ be some open and bounded domain. We consider a \emph{classical}\footnote{By \emph{classical} we mean defined in a pointwise manner.} linear operator $B_{\mu;\circ}$ on $\Omega$ with classical domain 
\begin{align*}
    D(B_{\mu;\circ}) =\{ x\in C({\xoverline\Omega}): x_{|\partial\Omega}=0, B_{\mu;\circ} x\in C(\Omega)\} 
\end{align*}
and aim at solving 
\begin{gather}
    B_{\mu;\circ} x_\mu = f_\mu
    \text{ (pointwise) on }\Omega, 
    \qquad
    {x_\mu}_{|\partial\Omega}=0.
    \label{eq:klassic_form}
\end{gather}
Note that the definition of $B_{\mu;\circ}$ also incorporates essential homogeneous boundary conditions (in case of a PDAE described below this is the initial condition, which is independent of the parameter). Let $\{B^*_{\mu;\circ},D(B^*_{\mu;\circ})\}$ denote the operator, which is adjoint to  $\{B_{\mu;\circ},D(B_{\mu;\circ})\}$, i.e., $B^*_{\mu;\circ}$ is defined as the formal adjoint  of $B_{\mu;\circ}$ by $\inner{B_{\mu;\circ} x, y}_{L_2(\Omega)} = \inner{x, B^*_{\mu;\circ} y}_{L_2(\Omega)}$ for all $x,y\in C^\infty_0(\Omega)$ and its domain $D(B^*_{\mu;\circ})$ which includes the corresponding adjoint essential boundary conditions (so that the above equation still holds true for all $x\in D(B_{\mu;\circ})$, $y\in D(B^*_{\mu;\circ})$). Denoting the range of an operator $B$ by $R(B)$, we have $B_{\mu;\circ}: D(B_{\mu;\circ}) \ra R(B_{\mu;\circ})$ and $B^*_{\mu;\circ}: D(B^*_{\mu;\circ}) \ra R(B^*_{\mu;\circ})$. The following assumptions\footnote{The framework in \cite{Dahmen2011} is slightly more general.} turned out to be crucial for ensuring the well-posedness:
\begin{compactenum}[\hskip8pt(B1)]
    \item $D(B_{\mu;\circ}), D(B^*_{\mu;\circ}), R(B^*_{\mu;\circ}) \subseteq  L_2(\Omega)$ with all embeddings being dense; \label{it:annahme1}
    \item $B^*_{\mu;\circ}$ is injective on $D(B^*_{\mu;\circ})$.\label{it:annahme2}
\end{compactenum}
Due to (B2), the injectivity of the adjoint operator, the following quantity 
\begin{gather*}
    \opnorm{\cdot}_\mu := \norm{B_{\mu}^* \cdot}_{L_2(\Omega)}
\end{gather*}
is a norm on $D(B_{\mu;\circ}^*)$, where $B_\mu^*$ is to be understood as the continuous extension of $B^*_{\mu;\circ}$ onto $Y_\mu$, i.e., $B_\mu^*: Y_\mu \ra L_2(\Omega)$, where
\begin{gather*}
    Y_\mu := \closure_{\opnorm{\cdot}_\mu}(D(B_{\mu;\circ}^*)),\quad 
    \inner{v,w}_{Y_\mu} := \inner{B_\mu^*v, B_\mu^*w}_{L_2(\Omega)},\quad
    \norm{v}_{Y_\mu}^2 := \inner{v,v}_{Y_\mu}=\opnorm{v}_\mu^2,
\end{gather*}
is a Hilbert space. Defining the bilinear form 
\begin{gather*}
    b_\mu: L_2(\Omega)\times Y_\mu \ra \R
    \quad\text{by}\quad
    b_\mu(x,y) := \inner{x,B_\mu^*y}_{L_2(\Omega)},
\end{gather*}
yields an ultraweak form of \eqref{eq:klassic_form}:
For $f\in Y_\mu'$\footnote{$Y_\mu'$ denotes the dual space of $Y_\mu$ w.r.t.\ the pivot space $L_2(\Omega)$.}, determine $x_\mu\in L_2(\Omega)$ such that
\begin{gather}
    b_\mu(x_\mu,y_\mu) = f_\mu(y_\mu)\qquad\forall y_\mu\in Y_\mu. \label{eq:weak_form}
\end{gather}
Well-posedness including optimal stability is now ensured:

\begin{lemma} \label{co:dahmen}
Problem \crefl{eq:weak_form} has a unique solution $x_\mu\in L_2(\Omega)$ and is optimally stable, i.e., $\gamma_\mu=\beta_\mu=\beta^*_\mu = 1$, where the continuity constant is defined as
\begin{gather*}
    \gamma_\mu := \sup_{x\in L_2(\Omega)} \sup_{y_\mu\in Y_\mu} \frac{b_\mu(x,y_\mu)}{\norm{x}_{L_2(\Omega)} \norm{y_\mu}_{Y_\mu}},
\end{gather*}
    and primal resp.\ dual inf-sup constants read
\begin{gather*}
    \beta_\mu:= \inf_{x\in L_2(\Omega)} \sup_{y_\mu\in Y_\mu} \frac{b_\mu(x,y_\mu)}{\norm{x}_{L_2(\Omega)} \norm{y_\mu}_{Y_\mu}},
    \qquad
    \beta^*_\mu := \inf_{y_\mu\in Y_\mu} \sup_{x\in L_2(\Omega)} \frac{b_\mu(x,y_\mu)}{\norm{x}_{L_2(\Omega)} \norm{y_\mu}_{Y_\mu}}.
\end{gather*}
\end{lemma}
\begin{proof}
See \cite[Proposition 3.1 and Corollary 3.2]{Dahmen2011}. 
\end{proof}

\subsection{An ultraweak formulation of PDAEs}\label{Sec:vwPDAE} 
We are now going to apply the framework of \S\ref{ch:weak} to the classical form \eqref{eq:PDAE} of the PDAE. Again, w.l.o.g. we restrict ourselves to homogeneous initial conditions $x_\mu(0) = 0$, as stated in \S\ref{Sec:Hom}. 

It is immediate that we can generalize ultraweak formulations for scalar-valued functions in $L_2(\Omega)$ as above to systems, i.e., $L_2(\Omega)^n \equiv L_2(\Omega;\R^n)$. For PDAEs, we choose $L_2(I)^n$ with the inner product $\inner{\cdot,\cdot}_{L_2} \equiv \inner{\cdot,\cdot}_{L_2(I)^n}$, whereas $(\cdot,\cdot)$ denotes the Euclidean inner product of vectors. The linear operator $\{B_{\mu;\circ}, D(B_{\mu;\circ})\}$ corresponding to \eqref{eq:PDAE} reads 
\begin{gather*}
    B_{\mu;\circ} := E \ts{\frac{d}{dt}} - A_\mu,
    \qquad
    D(B_{\mu;\circ}) := \braces{x \in C^{k_\mu}(I)^n \cap C(\bar I)^n: x(0) = 0}.
\end{gather*}
The formal adjoint operator $B_{\mu;\circ}^*$ is easily derived by integration by parts, i.e., 
\begin{align*}
    \inner{B_{\mu;\circ} x,y}_{L_2}
    &= \inner{E\dot x - A_\mu x,y}_{L_2}
     = \inner{\dot x,E^Ty}_{L_2} - \inner{x,A_\mu^Ty}_{L_2}\\
    &= \inner{x(T),E^T y(T)} - \inner{x(0), E^T y(0)} - \inner{x,E^T\dot y}_{L_2} - \inner{x,A_\mu^Ty}_{L_2}\\
    &= \inner{x,-E^T\dot y - A_\mu^Ty}_{L_2}
     =: \inner{x, B_{\mu;\circ}^* y}_{L_2} \quad \forall x,y\in C^\infty_0(I)^n,
\end{align*}
which shows that
\begin{gather}
    B_{\mu;\circ}^* := -E^T \ts{\frac{d}{dt}} - A_\mu^T,
        \nonumber\\
    D(B^*_{\mu;\circ}) \equiv C^1_{E}(I)^n := \braces{y \in C^1(I)^n  \cap C(\bar I)^n: y(T) \in \kernel(E^T)}.\label{eq:defC1E}
\end{gather}
In fact, $\inner{B_{\mu;\circ} x,y}_{L_2} = \inner{x, B^*_{\mu;\circ} y}_{L_2}$ for all $x \in D(B_{\mu;\circ})$ and $y\in D(B^*_{\mu;\circ})$ since the boundary terms above still vanish thanks to $x(0) = 0$ and $y(T) \in \kernel(E^T)$. Moreover
\begin{gather*}
    R(B_{\mu;\circ}) = C^{k_\mu-1}(I)^n \cap C(\bar I)^n,
    \qquad
    R(B^*_{\mu;\circ}) = C(\bar I)^n.
\end{gather*}

\begin{lemma}
    We have $D(B_{\mu;\circ}),D(B^*_{\mu;\circ}),R(B^*_{\mu;\circ}) \subset L_2(I)^n$ with dense embeddings.
\end{lemma}
\begin{proof}
By the definition of $H^1_0(I)^n$ and \cite[Cor.\ 7.24]{ArendtUrban} (for $H^1(I)^n$ instead of $H^1(I)$ there, which is a trivial extension), we have 
\begin{gather*}
    C^\infty_0(I)^n \subset H^1_0(I)^n \subset C(\bar{I})^n,
    \quad\text{hence}\quad
    C^\infty_0(I)^n = C^\infty_0(I)^n \cap C(\bar{I})^n.
\end{gather*}
With that, $C^\infty_0(I)^n \subseteq D(B_{\mu;\circ}),D(B^*_{\mu;\circ}),R(B^*_\circ) \subset L_2(I)^n$ is easy to see. Since $C^\infty_0(I)^n$ is dense in $L_2(I)^n$, its supersets $D(B_{\mu;\circ}),D(B^*_{\mu;\circ}),R(B^*_{\mu;\circ})$ are also dense in $L_2(I)^n$.
\end{proof}

The above lemma ensures assumption (B1). Next, we consider (B2).

\begin{lemma} \label{ch:injektivitaet}
The adjoint operator $\{B^*_{\mu;\circ}, D(B^*_{\mu;\circ})\}$ is injective, i.e., for $y_\mu,z_\mu\in D(B^*_{\mu;\circ})$ with $B_{\mu;\circ}^* y_\mu = B_{\mu;\circ}^* z_\mu$ we have $y_\mu = z_\mu$.
\end{lemma}
\begin{proof}
Setting $d_\mu := y_\mu - z_\mu$, we get $B^*_{\mu;\circ} d_\mu = 0$ and
\begin{align*}
    -E^T \dot d_\mu(t) - A_\mu^T d_\mu(t) = 0,\quad\forall t\in I,
    &\qquad
        d_\mu(T) = y_\mu(T) - z_\mu(T) \in \kernel(E^T).
\end{align*}
Due to regularity of $\braces{E,-A_\mu}$ (and thus also of $\braces{-E^T,-A_\mu^T}$), there are regular matrices $P_\mu$, $Q_\mu\in\C^{n\times n}$, which allow us to transform the problem into Weierstrass-Kronecker normal form, \cite{HairerWanner2,KunkelMehrmann}, i.e., 
\begin{gather*}
    P_\mu E^TQ_\mu = \begin{pmatrix} I\!d_m & 0 \\ 0 & N_\mu\end{pmatrix},
    \quad 
    P_\mu A_\mu^TQ_\mu = \begin{pmatrix} R_\mu & 0 \\ 0 & I\!d_{n-m}\end{pmatrix},
    \quad
    Q_\mu^{-1}d_\mu(t) = \begin{pmatrix} u_\mu(t) \\ v_\mu(t) \end{pmatrix},
\end{gather*}
where $I\!d_n\in\R^{n\times n}$ is the identity and $N_\mu$ is a nilpotent matrix with nilpotency index $k_\mu$. This yields the equivalent representation 
\begin{subequations}
\begin{align}
    \dot u_\mu(t) + R_\mu u_\mu(t) &= 0,\qquad\forall t\in I, \label{eq:wkn_ode}\\
    N_\mu \dot v_\mu(t) + v_\mu(t) &= 0,\qquad\forall t\in I,\label{eq:wkn_alg}\\
    Q_\mu\begin{pmatrix} u_\mu(T) \\ v_\mu(T) \end{pmatrix} 
    &\in \kernel(E^T) \label{eq:wkn_ew}.
\end{align}
\end{subequations}
The ODE \eqref{eq:wkn_ode} has the general solution  $u_\mu(t) = u_\mu(T)\, e^{-R_\mu(T-t)}$. By \eqref{eq:wkn_ew} we get  
\begin{gather*}
    E^T Q_\mu \begin{pmatrix} u_\mu(T) \\ v_\mu(T) \end{pmatrix} 
    = 0
    = P_\mu E^T Q_\mu \begin{pmatrix} u_\mu(T) \\ v_\mu(T) \end{pmatrix} 
    = \begin{pmatrix} I\!d_m & 0 \\ 0 & N_\mu\end{pmatrix} 
    \begin{pmatrix} u_\mu(T) \\ v_\mu(T) \end{pmatrix}
    = \begin{pmatrix} u_\mu(T) \\ N_\mu v_\mu(T) \end{pmatrix},
\end{gather*}
so that $u_\mu(T)=0$ and hence $u_\mu(t) = u_\mu(T)\, e^{-R_\mu(T-t)}=0$ for all $t\in I$.

The initial value problem $N_\mu \dot v_\mu(t) + v_\mu(t) = q_\mu(t)$, $t\in I$, $v_\mu(T)=v_{\mu,T}$ with some $q_\mu \in C^{k_\mu-1}(\bar I)^{n-m}$ has the unique solution $v_\mu(t) = \sum_{i=0}^{k_\mu-1} (-1)^iN_\mu^iq_\mu^{(i)}$, if the initial value $v_{\mu,T}$ is consistent, see e.g.\ \cite{brenan1995numerical}. We apply this for $q_\mu\equiv 0 \in C^{k_\mu-1}(\bar I)^{n-m}$. Then, by the solution formula, we get $v_\mu\equiv 0$, since the initial value in \crefl{eq:wkn_ew} is by definition trivially consistent. This yields $d_\mu \equiv 0$, i.e., $y_\mu = z_\mu$.
\end{proof}

Hence, we set $\opnorm{\cdot}_\mu := \norm{B_{\mu}^* \cdot}_{L_2}$ and choose trial and test spaces as 
\begin{gather}
    \label{eq:veryweakform}
    X := L_2(I)^n,
    \quad
    Y_\mu := \closure_{\opnorm{\cdot}_\mu}\parens[\big]{C^1_{E}(I)^n},
    \quad 
    b_\mu(x,y) := \inner{x, B_{\mu}^* y}_{L_2},
\end{gather}
see \eqref{eq:defC1E} and obtain the following result.

\begin{lemma}\label{Lem:Y}
    Under the above assumptions, we have for all $\mu\in\cP$ that $Y_\mu\equiv Y$, where
    \begin{align*}
        Y:= H^1_E(I)^n:=\{ v\in H^1(I)^n:\, v(T)\in\kernel(E^T)\}.
    \end{align*}
\end{lemma}
\begin{proof}
    Clearly $C^1_{E}(I)^n\subset H^1_E(I)^n$, so that $Y_\mu\subseteq Y$ for all $\mu\in\cP$. Now, let $y\in Y = H^1_E(I)^n$, then, by density, there is a sequence $(y_\ell)_{\ell\in\N}\subset C^1_{E}(I)^n$ such that $\| y_\ell-y\|_{H^1(I)^n}\to 0$ as $\ell\to 0$. Since $\cP$ is compact, we have that
    \begin{align*}
        \opnorm{y_\ell-y}_\mu
        &= \| E^T (\dot{y}_\ell-\dot{y}) 
        + A_\mu^T(y_\ell-y)\|_{L_2}
        \le \max\{ \| E\|, \|A_\mu\|\}\, \| y_\ell-y\|_{H^1(I)^n}
        \to 0
    \end{align*}
    as $\ell\to\infty$. Hence, $y\in \closure_{\opnorm{\cdot}_\mu}\parens[\big]{C^1_{E}(I)^n}=Y_\mu$, i.e. $Y \subseteq Y_\mu$.
\end{proof}

The latter result must be properly interpreted. It says that $Y_\mu$ and $Y$ coincide as sets. However, the norm $\opnorm{\cdot}_\mu$ (and thus the topology) still depends on the parameter. The same holds true for the dual space $Y'$ of $Y$ induced by the $L_2$-inner product and normed by
\begin{align*}
    \opnorm{f}_{\mu}':= \sup_{y\in Y} 
         \frac{(f,y)_{L_2}}{\opnorm{y}_\mu}.
\end{align*}
In particular, we have a generalized Cauchy-Schwarz inequality $(f,y)_{L_2}\le \opnorm{f}_{\mu}'\, \opnorm{y}_{\mu}$.

\begin{lemma}\label{lm:wp_vwDAE}
    Let $f_\mu\in Y'$. Then, there exists a unique weak solution $x_\mu\in X$ of 
    \begin{gather}
        \label{eq:vwDAE}
        b_\mu(x_\mu,y)=f_\mu(y),\qquad\forall y\in Y.
    \end{gather}
    If  \eqref{eq:PDAE} admits a classical solution, then it coincides with $x_\mu$. 
    Moreover, $\gamma_\mu=\beta_\mu=\beta_\mu^* = 1$ for the constants defined in Lemma \ref{co:dahmen}. 
\end{lemma}
\begin{proof}
The existence of a unique solution $x_\mu\in X$ (as well as $\gamma_\mu=\beta_\mu=\beta_\mu^* = 1$) is an immediate consequence of Lemma \ref{co:dahmen}. It only remains to show that $x_\mu$ satisfying \eqref{eq:vwDAE} is a weak solution of \eqref{eq:PDAE}. To this end, let $\tilde{f}_\mu \in C(I)^n$ be given such that there exists a classical solution $\tilde{x}_\mu\in C^1(\bar I)^n$ with $B_{\mu;\circ} \tilde{x}_\mu(t) = \tilde{f}_\mu(t),\forall t\in I$ and $\tilde{x}_\mu(0) = 0$. Then, define $f_\mu\in Y'$ by $f_\mu(y)  := \inner{\tilde{f}_\mu,y}_{L_2}$. We need to show that the classical solution $\tilde{x}_\mu$ of \eqref{eq:PDAE} is also the unique solution of \eqref{eq:vwDAE}. First, for $y\in C^1_E(I)^n$, integration by parts yields 
$b_\mu(\tilde{x}_\mu,y) - f_\mu(y)
    = \inner{\tilde{x}_\mu,B_\mu^*y}_{L_2} - f_\mu(y)
    = \inner{B_{\mu;\circ} \tilde{x}_\mu -\tilde{f}_\mu,y}_{L_2}
    = 0$. 
Second, let $y\in Y\setminus C^1_{E}(I)^n$, then there is  $(\tilde y_{\ell})_{{\ell}\in\N}\subset C^1_{E}(I)^n$  converging to $y$ in $Y$, i.e.,  $\lim_{{\ell}\ra\infty} \opnorm{y - \tilde y_{\ell}}_{\mu} = 0$. Then, by the generalized Cauchy-Schwarz inequality
\begin{align*}
    \abs{b_\mu(\tilde{x}_\mu,y) - f_\mu(y)}
    &= \abs{b_\mu(\tilde{x}_\mu,y) - f_\mu(y) - b_\mu(\tilde{x}_\mu,\tilde y_{\ell}) + f_\mu(\tilde y_{\ell})} \\
    &= \abs{\inner{\tilde{x}_\mu,B_\mu^*(y-\tilde y_{\ell})}_{L_2} - f_\mu(y-\tilde y_{\ell})}\\
    &\leq \norm{\tilde{x}_\mu}_{L_2} \norm{B_\mu^*(y-\tilde y_{\ell})}_{L_2}
        + \opnorm{f_\mu}_{\mu}'\, \opnorm{y-\tilde y_{\ell}}_{\mu}\\
    &= \parens[\big]{\norm{\tilde{x}_\mu}_{L_2} + \opnorm{f_\mu}_{\mu}'}\, \opnorm{y-\tilde y_{\ell}}_{\mu} \ra 0
    \quad \text{ as } {\ell}\to\infty,
\end{align*}
so that \eqref{eq:vwDAE} holds for $\tilde{x}_\mu$. 
\end{proof}

For the ultraweak PDAE \eqref{eq:vwDAE}, we need a right-hand side $f_\mu\in Y'$. However, typically, the right-hand side is given within context of \eqref{eq:PDAE} as a function of time, i.e., $g_\mu: I \ra \R^n$. Then, we simply define $f_\mu\in Y'$ by
\begin{align}
    f_\mu(y) := \inner{g_\mu,y}_{L_2} = \bintegral{I}{}{\inner{g_\mu(t),y(t)}}{t},
    \quad y\in Y.
    \label{eq:construct_rhs}
\end{align}

\section{Petrov-Galerkin discretization}\label{Sec:3}

The next step towards a numerical method for solving an ultraweak operator equation is to introduce finite-dimensional trial and test spaces yielding a Petrov-Galerkin discretization. In this section, we shall first review Petrov-Galerkin methods in general terms and then detail the specification for PDAEs.

\subsection{Petrov-Galerkin method}
\label{Subsec:3.1}
In order to determine a numerical approximation, we are going to construct an appropriate finite-dimensional trial space $X_\mu^\cN\subset X=L_2(I)^n$ and a parameter-independent test space $Y^\cN \subset Y$  of finite (but possibly large) dimension $\cN\in\N$. Then, we are seeking $x_\mu^\cN\in X_\mu^\cN$ such that
\begin{gather}
    \label{eq:vwDAE_disc}
    \quad b_\mu(x_\mu^\cN,y^\cN) = f_\mu(y^\cN), \qquad \forall y^\cN\in Y^\cN,
\end{gather}
which leads to solving a linear system of equations $\bm{B}_\mu^\cN\bm{x}_\mu^\cN = \bm{f}_\mu^\cN$ in $\R^\cN$.

\begin{remark}\label{Rem:1}
    \begin{compactenum}[(a)]
        \item If one would choose a discretization with $\dim(Y^\cN)>\dim(X_\mu^\cN)$, one would need to solve a least squares problem $\|\bm{B}_\mu^\cN\bm{x}_\mu^\cN - \bm{f}_\mu^\cN\|^2\to\min$.
        \item\label{it:discret_constants} If one defines the trial space according to $X_\mu^\cN := B_\mu^* Y^\cN$, then it is easily seen that the discrete problem \crefl{eq:vwDAE_disc} is well-posed and optimally conditioned, \cite{BSU19}, i.e., 
        \begin{align*}
            \gamma_\mu^\cN &:= \sup_{x\in X_\mu^\cN} \sup_{y\in Y^\cN} 
                \frac{b_\mu(x,y)}{\norm{x}_{L_2}\, \opnorm{y}_{\mu}}=1,\\
            \beta_\mu^\cN &:= \inf_{x\in X_\mu^\cN} \sup_{y\in Y^\cN} 
                \frac{b_\mu(x,y)}{\norm{x}_{L_2} \,\opnorm{y}_{\mu}}=1,
                \qquad 
            \beta_\mu^{*,\cN} := \inf_{y\in Y^\cN} \sup_{x\in X_\mu^\cN} 
                \frac{b_\mu(x,y)}{\norm{x}_{L_2}\, \opnorm{y}_{\mu}}=1.
        \end{align*}
        \item The Xu-Zikatanov lemma (\cite{Xu2003}) ensures that the Petrov-Galerkin error is comparable with the error of the best approximation, namely 
        \begin{gather}
            \label{eq:XuZik}
            \norm{x_\mu - x_\mu^\cN}_{L_2} \leq \frac{\gamma_\mu}{\beta_\mu^\cN} 
                \inf_{v^\cN\in X_\mu^\cN} \norm{x_\mu - v^\cN}_{L_2},
        \end{gather}
        so that the Petrov-Galerkin approximation is the best approximation (i.e., an \emph{identity}) for $\gamma_\mu=\beta_\mu^\cN=1$.
    \end{compactenum}
\end{remark}

The Petrov-Galerkin framework induces a residual-based error estimation in a straightforward manner. To describe it, let us recall that the \emph{residual} is defined for some $\tilde x\in L_2(I)^n$ as
\begin{gather*}
    r(\tilde x)\in Y',
    \qquad  r(\tilde x)[y] := 
        f_\mu(y) - b_\mu(\tilde{x},y),
    \quad y\in Y.
\end{gather*}
Then, it is a standard estimate that
\begin{align}\label{eq:ErrEst}
    \| x_\mu - x_\mu^\cN\|_{L_2} 
    &\le \frac{1}{\beta_\mu} \sup_{y\in Y} \frac{b_\mu(x_\mu - x_\mu^\cN, y)}{\opnorm{y}_{\mu}}
        = \frac{1}{\beta_\mu} \opnorm{r(x_\mu^\cN)}_{\mu}'
        =: \Delta_\mu^\cN
\end{align}
and $\Delta_\mu^\cN$ is a residual-based error estimator. Note that for $\beta_\mu = 1$ we have a error-residual identity $\| x_\mu - x_\mu^\cN\|_{L_2} = \opnorm{ r(x_\mu^\cN)}_{\mu}' = \Delta_\mu^\cN$.

\subsection{PDAE Petrov-Galerkin Discretization} \label{Sec:Disc}
\begin{figure}
\centering
\begin{tikzpicture}
       \tikzmath{
       		\h = 1.4; 
       		\n = 8; 
       		\dt = 1.4; 
       		\b = \n*\dt;
       }%
       \tikzstyle{myFunBesch}=[anchor=south west]
       \tikzstyle{myFun}=[thick]
       \tikzstyle{myAxis}=[]
       \draw[->,myAxis] (\b, 0)--(\b+0.5, 0) node[anchor=north] {$t$};
       \draw[->,myAxis] (0, 0)--(0, \h+0.5) node[anchor=north] {};
       \draw (1pt,\h)--(-3pt,\h) node[anchor=east] {$1$};
       \foreach \x/\m/\c/\l/\ll/\t/\f/\a in 
       	{ 0/1/black/-/-/{$0$}/{$\sigma_0$}/-,
       	  1/2/black/dashed/dashed/{$t_1$}/{}/dashed,
       	  2/2/black/dashed/dashed/{$t_{k-2}$}/{}/-,
       	  3/2/black/-/-/{$t_{k-1}$}/{$\sigma_{k-1}$}/-,
       	  4/2/black/-/-/{$t_k$}/{$\sigma_{k}$}/-,
       	  5/2/black/-/-/{$t_{k+1}$}/{$\sigma_{k+1}$}/-,
       	  6/2/black/dashed/dashed/{$t_{k+2}$}/{}/dashed,
       	  7/2/black/dashed/dashed/{$t_{K-1}$}/{}/-,
       	 \n/3/black/-/-/{$T$}/{$\sigma_K$}/-%
       	 }{%
       	 	\draw (\x*\dt,1pt)--(\x*\dt,-3pt) node[anchor=north] {\t};     	 	
     		\tikzmath{%
     			if \x < \n then {%
     				{\draw[\a,myAxis] (\x*\dt,0)--(\x*\dt+\dt,0);};%
     			} else {};%
     			if \m > 0 then {%
       			if \m == 1 then {%
       				{\draw[\c,\l,myFun] (\x*\dt+\dt,0)--(\x*\dt,\h) node[myFunBesch] {\f};};%
       			} else {%
       				if \m == 2 then {%
       					{\draw[\c,\l,myFun] (\x*\dt+\dt,0)--(\x*\dt,\h) node[myFunBesch] {\f};%
       					 \draw[\c,\ll,myFun] (\x*\dt-\dt,0)--(\x*\dt,\h);};%
       				} else {%
       					{\draw[\c,\ll,myFun] (\x*\dt-\dt,0)--(\x*\dt,\h) node[myFunBesch] {\f};};%
       				};%
       			};%
       			} else {};%
       		}%
       }
\end{tikzpicture}
\caption{Piecewise linear temporal discretization (hat functions).}\label{fig:hatfun}
\end{figure}
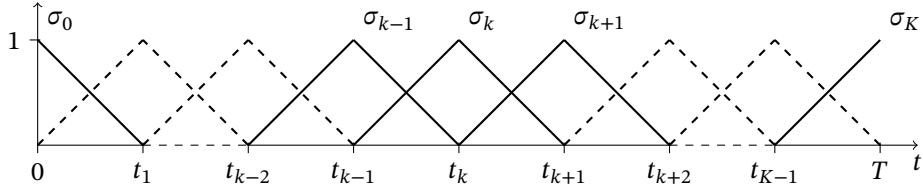

We are now going to specify the above general framework to PDAEs. This means that we need to introduce a suitable discretization in time. We fix a constant time step size $\Delta t := T/K $ (i.e., $K\in\N$ is the number of time intervals) and choose for simplicity equidistant nodes $t_k := k\Delta t,\, k=0,...,K$ in $I$. Denote by $\sigma_k$, $k=0,...,K$ piecewise linear splines corresponding to the nodes $t_{k-1}$, $t_k$ and $t_{k+1}$, see Figure \ref{fig:hatfun}. For $k\in\{ 0,K\}$, the hat functions are restricted to the interval $\bar{I}$. For realizing a discretization of higher order, one could simply use splines of higher degree.

As in \cite{BSU19}, we start by defining the test space and then construct inf-sup optimal trial spaces. To this end, let  $d := \dim (\kernel E^T)$ and assume that we have a basis $\braces{v_1,\dots,v_{d}}$ of $\kernel E^T$ at hand\footnote{This is in fact the reason why we restricted ourselves to parameter-independent matrices $E$ instead of $E_\mu$. We would then need to have a parameter-dependent basis for $\kernel E_\mu^T$, which is of course possible, but causes a quite heavy notation.}  and form a matrix $V:=(v_1,...,v_d)\in\R^{n\times d}$ by arranging the vectors as columns of $V$. 
Then, we construct $Y^\cN\subset Y=H^1_E(I)^n$ independent of the parameter and choose the trial space as $X_\mu^\cN:=B^*_\mu Y^\cN$, which will then guarantee that $\beta_\mu^\cN=1$. We suggest a piecewise linear discretization by
\begin{align*}
    Y^\cN :=\spann\braces[\big]{e_i \sigma_k:\ k=0,...,K-1,\; i=1,...,n} 
    \,\oplus\, 
    \spann\braces[\big]{v_i \sigma_K: \, i=1,\dots,d} \subset Y,
\end{align*}
where $e_i\in\R^n$ denotes $i$-th canonical vector. Then, we set
\begin{align*}
    X^\cN_\mu := B^*_\mu Y^\cN=&\spann\braces[\big]{-E^Te_i \dot\sigma_k -A_\mu^Te_i \sigma_k:\, k=0,...,K-1,\; i=1,...,n}\\
    & \oplus \spann\braces[\big]{-E^Tv_i \dot\sigma_K -A_\mu^Tv_i \sigma_K: \; i=1,\dots,d} \subset X=L_2(I)^n,
\end{align*}
with dimensions $\cN := \dim(X^\cN) = \dim(Y^\cN) = nK + d$. Then, Lemma \ref{lm:wp_vwDAE} and Remark \ref{Rem:1} ensure $\beta_\mu^\cN=\beta_\mu=\gamma_\mu=1$ and thus
\begin{align}
    \| x_\mu - x_\mu^\cN\|_{L_2} = \inf_{v^\cN\in X_\mu^\cN} \norm{x_\mu - v^\cN}_{L_2} = \opnorm{ r(x_\mu^\cN)}_{\mu}' = \Delta_\mu^\cN. \label{eq:res_identity}
\end{align}

\subsubsection{The linear system}
To construct the discrete linear system $\bm{B}_\mu^\cN\bm{x}_\mu^\cN = \bm{f}_\mu^\cN$ we need bases $\braces{\xi_1(\mu),\dots,\xi_\cN(\mu)}$ of $X^\cN_\mu$ and $\braces{\psi_1,\dots,\psi_\cN}$ of $Y^\cN$. 
The stiffness matrix ${\bm{B}}_\mu^\cN\in\R^{\cN\times \cN}$ can be computed by $[{\bm{B}}_\mu^\cN]_{j,i} := b_\mu(\xi_i(\mu), \psi_j) = (\xi_i(\mu), B_\mu^* \psi_j)_{L_2}$. We recall that $X^\cN_\mu = B^*_\mu Y^\cN$, which implies that $\xi_i(\mu)=B^*_\mu \psi_i$, so that $[{\bm{B}}_\mu^\cN]_{j,i}= (B^*_\mu \psi_i, B_\mu^* \psi_j)_{L_2}$ and ${\bm{B}}_\mu^\cN$ is in fact  symmetric positive definite. 

The right-hand side $\bm{f}_\mu^\cN \in\R^\cN$ reads  $\brackets{\bm{f}_\mu^\cN}_{j} := f_\mu(\psi_j)$. The discrete solution then reads $x_\mu^\cN := \sum_{i=1}^\cN \brackets{\bm{x}_\mu^\cN}_i\, \xi_i(\mu)$.

Recalling the finite element functions $\sigma_k$ in Figure \ref{fig:hatfun}, we define the inner product matrices for $k,\ell=0,...,K$ by
\begin{gather*}
\brackets{\bm{K}_{\Delta t}}_{k,\ell} := \inner{\dot\sigma_k,\dot\sigma_\ell}_{L_2(I)},\qquad
\brackets{\bm{L}_{\Delta t}}_{k,\ell} := \inner{\sigma_k,\sigma_\ell}_{L_2(I)},\qquad
\brackets{\bm{O}_{\Delta t}}_{k,\ell} := \inner{\dot\sigma_k,\sigma_\ell}_{L_2(I)},
\end{gather*}
and subdivide the matrices $\bm{\Pi}_{\Delta t} \in\R^{(K+1) \times (K+1)}$ for $\bm\Pi_{\Delta t}\in\braces{\bm{K}_{\Delta t},\bm{L}_{\Delta t},\bm{O}_{\Delta t}}$ according to
\begin{gather*}
\renewcommand*{\arraystretch}{1.5}
\bm\Pi_{\Delta t} 
= \begin{pmatrix}[cc]
    \bm\Pi_{\Delta t}^{1,1} & \bm\Pi_{\Delta t}^{1,2}\\ 
    \bm\Pi_{\Delta t}^{2,1} & \bm\Pi_{\Delta t}^{2,2}
    \end{pmatrix},\quad
\begin{split}
    \begin{alignedat}{4}
    &\bm\Pi_{\Delta t}^{1,1}  &&\in\R^{K \times K},\quad
    &&\bm\Pi_{\Delta t}^{1,2} &&\in\R^{K \times 1},\\
    &\bm\Pi_{\Delta t}^{2,1}  &&\in\R^{1 \times K},\quad
    &&\bm\Pi_{\Delta t}^{2,2} &&\in\R.
\end{alignedat}
\end{split}
\end{gather*} 
Then, the stiffness matrix also has block structure
\begin{gather*}
\renewcommand*{\arraystretch}{1.5}
    {\bm{B}}_\mu^\cN 
    = \begin{pmatrix}[cc]
    {\bm{B}}_\mu^{1,1} & {\bm{B}}_\mu^{1,2}\\ 
    {\bm{B}}_\mu^{2,1} & {\bm{B}}_\mu^{2,2}
    \end{pmatrix} 
    \in\R^{\mathcal{N}\times\mathcal{N}}
\end{gather*}
in form of Kronecker products of matrices, i.e. (with $V= (v_1,\dots,v_d) \in\R^{n\times d}$ as above), 
\begin{align*}
&{\bm{B}}_\mu^{1,1} = 
    \bm{K}_{\Delta t}^{1,1} \otimes EE^T
    +\bm{O}_{\Delta t}^{1,1} \otimes EA_\mu^T
    +(\bm{O}_{\Delta t}^{1,1})^T \otimes A_\mu E^T
    +\bm{L}_{\Delta t}^{1,1} \otimes A_\mu A_\mu^T 
    \in \R^{nK\times nK},\\
&{\bm{B}}_\mu^{1,2} = 
   \bm{O}_{\Delta t}^{1,2} \otimes EA_\mu^TV
   +\bm{L}_{\Delta t}^{1,2} \otimes A_\mu A_\mu^TV 
   \in\R^{nK\times d},\\
&{\bm{B}}_\mu^{2,1} = 
   (\bm{O}_{\Delta t}^{1,2})^T \otimes V^TA_\mu E^T
   +\bm{L}_{\Delta t}^{2,1} \otimes V^TA_\mu A_\mu^T \in\R^{d\times nK},\\
&{\bm{B}}_\mu^{2,2} 
    = \bm{L}_{\Delta t}^{2,2} \otimes V^TA_\mu A_\mu^TV 
    \in\R^{d\times d}.
\end{align*}

For the right-hand side, given some function $f_\mu : \bar I \ra \R^n$, we obtain a discretization $\bm{f}_\mu^\cN\in\R^{\mathcal{N}}$ in the sense of \eqref{eq:construct_rhs} by $\brackets{\bm{f}_\mu^\cN}_{j} =  \sum_{k=0}^K \inner{f_\mu(t_k)\, \sigma_k, \psi_j}_{L_2}$, $j=1,...,\mathcal{N}$. This means that we discretize $f_\mu$ in time by means of piecewise linears.  
Collecting the sample values of $f_\mu$ in one vector, i.e., $\bm{f}_{\mu,\Delta t} := \parens[\big]{{f}_\mu(t_0),\dots, {f}_\mu(t_K)}^T \in\R^{n(K+1)}$ we get that 
\begin{gather*}
\bm{f}_\mu^\cN = \bm{F}_{\Delta t}^T \bm{f}_{\mu,\Delta t}
\quad\text{where}\quad
\renewcommand*{\arraystretch}{1.5}
    \bm{F}_{\Delta t} := \begin{pmatrix}[cc]
    I\!d_n \otimes \bm{L}_{\Delta t}^{1,1} & V \otimes \bm{L}_{\Delta t}^{1,2} \\
    I\!d_n \otimes \bm{L}_{\Delta t}^{2,1} & V \otimes \bm{L}_{\Delta t}^{2,2}
    \end{pmatrix} \in \R^{n(K+1) \times \cN}
\end{gather*}
and $I\!d_n\in\R^{n\times n}$ again denoting the $n$-dimensional identity matrix. 

As already noted above, of course, one could use different discretizations (e.g.\ higher order or different discretizations for $f_\mu$ and the test functions) and we choose the described one just for simplicity.

The efficient numerical solution of this linear system requires a solver that takes the specific structure into account. For similar systems arising from space-time (ultraweak) variational formulations of heat, transport and wave equations, such specific efficient solvers have been introduced in \cite{HPSU20,Henning2021}. The structure of the system above is different and we will consider the development of efficient solvers in future research, see Section \ref{Sec:6}.

\subsubsection{Special case: Linear DAEs}\label{Sec:flDAE}
We are going to specify the above general setting to the special case of fully linear DAEs, namely 
\begin{gather}
    E\dot{x}(t) - Ax(t) = Bu(t)+g_{\mu_2}(t),\quad t\in I, 
    \qquad
    x(0)=0,
    \label{eq:flDEA}
\end{gather}
in which the right-hand side is given in terms of a matrix $B\in \R^{n\times m}$, a control $u:\bar{I} \ra \R^m$, $m$ denoting some input dimension and a function $g_{\mu_2}:\bar{I}\to\R^n$, which arises from the reduction to homogeneous initial conditions, see \S\ref{Sec:Hom}. The initial condition is assumed to be parameterized through $g_{\mu_2}$ by $\mu_2\in\cP_2\subset\R^{P_2}$, $P_2\in\N$. 
In view of \eqref{eq:affinedecomp} and \S\ref{Sec:Hom}, we get
\begin{align*}
	g_{\mu_2}(t)
	&= \sum_{q=1}^{Q_x} \vartheta^x_q(\mu_2)\, \big(A\bar{x}_q(t) - E\dot{\bar{x}}_q(t)\big)
	=: \sum_{q=1}^{Q_x} \vartheta^x_q(\mu_2)\, z_q(t),
\end{align*}
where $\bar{x}_q\in C^1(\bar{I})^n$ are smooth extensions of $\tilde x_{0,q}$, i.e., $\bar{x}_q(0) = \tilde x_{0,q}$, $q=1,...,Q_x$. 

We view the control and the initial condition (via $g_{\mu_2}$) as parameters, i.e., $f_\mu(t) = B{\mu_1}(t)+g_{\mu_2}(t)$, $\mu=(\mu_1,\mu_2)$, which means that the parameter set would be infinite-dimensional and needs to be discretized. 
Using the same kind of discretization as above, we can use the samples of the control as parameter, i.e., 
\begin{align*}
    \bm\mu_1 := \parens[\big]{u(t_0),...,u(t_K)}^T \in\cP_1=\R^{P_1},
    \quad
    P_1=m(K+1),
\end{align*}
and similar for the initial condition $\bm{z}_q := \parens[\big]{z_q(t_0),..., z_q(t_K)}^T \in\R^{n(K+1)}$, $q=1,...,Q_x$. Then, we get 
\begin{align*}
    \bm{f}_\mu^\cN=\bm{F}_{\Delta t}^T(B\otimes I\!d_{K+1})\, \bm\mu_1 
    + \sum_{q=1}^{Q_x} \vartheta_q^x(\mu_2)\, \bm{F}_{\Delta t}^T\bm{z}_q,
\end{align*}
so that the parameter dimension is $P=P_1+P_2=m(K+1)+P_2$, which might be large. The right-hand side $\bm{f}_\mu^\cN$ thus also admits an affine decomposition with $Q_f=P_1+Q_x=m(K+1)+Q_x$ terms. However, this number might be an issue concerning efficiency if $K$ is large. Nevertheless, if $m\ll n$, we still have $Q_f \ll \cN$.

Moreover, in the linear case, the matrix $A_\mu\equiv A$ is independent of the parameter, which means (among other facts) that trial and test spaces are parameter-independent, as sets and also w.r.t.\ their topology. Note that this is the most common case for system theoretic MOR methods (like BT), which are often even restricted to this case, \cite{MehrmannStykel2005,Stykel2003}, with the exception \cite{MR4171473}. Our setting seems more flexible in this regard and fully linear DAEs are just a special case.

\section{Model order reduction: The Reduced Basis Method}\label{Sec:MOR}
The Reduced Basis Method (RBM) is a model order reduction technique which has originally been constructed for parameterized partial differential equations (PPDEs), see e.g.\ \cite{BennerSIAM2017,rozza2016,quarteroni2015}. In an offline training phase, a reduced basis of size $N\ll\cN$ is constructed (typically in a greedy manner, see Algorithm \ref{Alg:Greedy} below) from sufficiently detailed approximations for certain parameter samples (also called \enquote{truth} approximations or \emph{snapshots}), which are computed e.g.\ by a Petrov-Galerkin method as described above. 
In particular, $\cN$ is assumed to be sufficiently large in order to ensure that $x_\mu^\cN$ is (at least numerically) indistinguishable from the exact state $x_\mu$, which explains the name \enquote{truth}. As long as an efficient solver for the detailed problem is available, we may assume that the snapshots can be computed in $\mathcal{O}(\cN)$ complexity.

Given some parameter value $\mu$, the reduced approximation $x_N(\mu)$\footnote{For all quantities of the reduced system, we write the parameter $\mu$ as an argument in order to clearly distinguish the detailed approximation $x_\mu^\cN$ from the reduced approximation $x_N(\mu)$ for the same parameter.}  is then computed by solving a reduced system of dimension $N$.  
Thanks to the affine decomposition \eqref{eq:affinedecomp}, several quantities for the reduced system can be precomputed and stored so that a reduced approximation is determined in $\mathcal{O}(N^3)$ operations, independent of $\cN$ (which is called \emph{online efficient}). Moreover, an \emph{a posteriori} error estimator $\Delta_N(\mu)$ guarantees a certification in terms of an online efficiently computable upper bound for the error, i.e., $\| x_\mu^\cN-x_N(\mu)\|_{L_2}\le \Delta_N(\mu)$.

We are going to use this framework for PDAEs of the form \eqref{eq:PDAE}. Model reduction of \eqref{eq:PDAE} may be concerned (at least) with the following quantities
\begin{compactitem}
    \item size $n$ of the system,
    \item dimension $K$ of the temporal discretization,
\end{compactitem}
where we have in mind to solve \eqref{eq:PDAE} extremely fast for several values of the parameter $\mu$. 
As mentioned earlier, the first issue has extensively been studied in the literature e.g.\ by system theoretic methods, in particular for fully linear DAEs \eqref{eq:flDEA}. This can be done independently from the subsequent reduction w.r.t.\ $K$ (both for parameterized and non-parameterized versions), so that we even might assume that such Model Order Reduction (MOR) techniques have already been applied in a preprocessing step. We mention  \cite{MR4171473} for a system theoretic MOR for parameter-dependent DAEs. Here, we are going to consider the reduction w.r.t.\ time using the RBM based upon a variational formulation w.r.t.\ the time variable.

We restrict ourselves to the reduction of the fully linear case of \eqref{eq:flDEA} as it easily shows how the RBM-inspired model reduction can be combined with existing system theoretic approaches to reduce the size of the system (e.g.\ in a preprocessing step). In the fully linear case, the matrix $A_\mu \equiv A$ and hence all operators and bilinear forms on the left-hand side are parameter-independent. This implies in addition that the ansatz space $X^\cN_\mu \equiv X^\cN$ and the norm $\opnorm{\cdot}_\mu \equiv \opnorm{\cdot}$ inducing the topology on the test space are parameter-independent as well, which of course simplifies the framework. However, parameter-dependent matrices $A_\mu$ can be treated similar to the RBM for ultraweak formulations of PPDEs as described e.g.\ in \cite{BSU19,Henning2021,UrbanCIME2022}. However, we note that the RB approach also allows the treatment of more general PDAEs and is not restricted to fully linear  systems \eqref{eq:flDEA}, in particular w.r.t.\ the right-hand side.

The idea of the RBM can be described as follows: One determines sample values 
\begin{align*}
    S_N:=\{\mu^{(1)},...,\mu^{(N)}\}\subset\cP
\end{align*} 
of the parameters in an offline training phase by a greedy procedure described in Algorithm \ref{Alg:Greedy} below. Then, for each $\mu\in S_N$, we determine a sufficiently detailed \enquote{snapshot} $x_\mu^\cN\in X^\cN$ by the ultraweak Petrov-Galerkin discretization as in \S\crefl{Sec:Disc} and obtain a reduced space of dimension $N$ by setting
\begin{align*}
    X_N := \spann\braces{x_\mu^\cN:\, \mu\in S_N}
    =: \spann\braces{\zeta_1,...,\zeta_N} 
    \subset X^\cN.
\end{align*}
We also need a reduced \emph{test} space for the Petrov-Galerkin method. Recalling that the operator is parameter-independent here ($B_\mu\equiv B$) and also the trial space $X^\cN$ is independent of $\mu$, we can easily identify the optimal test space. In fact, for each snapshot there exists a unique $y_\mu^\cN\in Y^\cN$ such that $x_\mu^\cN=B^* y_\mu^\cN$. Then, we define
\begin{align*}
    Y_N := \spann\braces{y_\mu^\cN:\, \mu\in S_N}
        =: \spann\braces{\eta_1,...,\eta_N}
        \subset Y^\cN.\footnotemark
\end{align*}
\footnotetext{For efficiency reasons, in fact, we first determine $\eta_i$ and then simply set $\zeta_i:=B^*\eta_i$.}
Then, given a new parameter value $\mu\in\cP$, one determines the \emph{reduced approximation} $x_N(\mu)\in X_N$ by solving (recall that here $b_\mu\equiv b$)
\begin{align*}
    b(x_N(\mu), y_N) = f_\mu(y_N)
    \quad\text{ for all } y_N\in Y_N.
\end{align*}
If $N\ll \cN = nK+d$, we can compute a reduced approximation with significantly less effort as compared to the Petrov-Galerkin (or a time-stepping) method. To determine the reduced approximation $x_N(\mu)$, we have to solve a linear system of the form $\bm{B}_N \bm{x}_N(\mu)=\bm{f}_N(\mu)$, where the stiffness matrix is given by $[\bm{B}_N]_{j,i} = b(\zeta_i,\eta_j)$, $i,j=1,...,N$, recalling that the bilinear form is parameter-independent. Hence, $\bm{B}_N\in\R^{N\times N}$ can be computed and stored in the offline phase. 
For the right-hand side, we use the affine decomposition \eqref{eq:affinedecomp} and get
$[\bm{f}_N(\mu)]_j = \sum_{q=1}^{Q_f} \vartheta_g^f(\mu)\, (\tilde{f}_q, \eta_j)_{L_2}$.
The quantities $(\tilde{f}_q, \eta_j)_{L_2}$ can be precomputed and stored in the offline phase, so that $\bm{f}_N(\mu)$ is computed online efficient in $\mathcal{O}(Q_fN)$ operations.  
Obtaining the coefficient vector $\bm{x}_N(\mu)$, the reduced approximation results in $x_N(\mu) = \sum_{i=1}^N [\bm{x}_N(\mu)]_i\, \zeta_i$.  
Note that the matrix $\bm{B}_N$ is typically densely populated so that the numerical solution requires in general $\mathcal{O}(N^3)$ operations. 

The announced greedy selection of the samples is based upon the residual error estimate (here identity) $\Delta_\mu^\cN$ in \eqref{eq:ErrEst} resp.\ \eqref{eq:res_identity} for the reduced system described as follows: 
In a similar manner as deriving $\Delta_\mu^\cN$ in \eqref{eq:ErrEst} we get a residual based error estimator for the reduced approximation 
\begin{align*}
    \| x_\mu^\cN - x_N(\mu)\|_{L_2} 
    &\le \frac{1}{\beta^\cN} \sup_{y\in Y^\cN} 
    \frac{b(x_\mu^\cN - x_N(\mu), y)}{\opnorm{y}}
        = \frac{1}{\beta^\cN} \opnorm{r(x_N(\mu))}'
        =: \Delta_N(\mu),
\end{align*}
since the bilinear form and the norm in $Y$ are parameter-independent here.
Hence, the inf-sup constant is parameter-independent as well, i.e., $\beta_\mu^\cN \equiv \beta^\cN$ and it is unity by Remark \ref{Rem:1}, so that 
\begin{align}\label{eq:ErrEstRB}
	\| x_\mu^\cN - x_N(\mu)\|_{L_2} = \opnorm{r(x_N(\mu))}' = \Delta_N(\mu). 
\end{align}
Its computation can be done in an online efficient manner in $\mathcal{O}(N)$ operations by determining Riesz representations in the offline phase, see  \cite{BennerSIAM2017,rozza2016,quarteroni2015}. 
We use this error identity in the greedy method in Algorithm \ref{Alg:Greedy} below.

\begin{algorithm}
    \caption{(Weak) Greedy method}\label{Alg:Greedy}%
    \begin{algorithmic}[1] 
	\Statex \textbf{input:} training sample $\cP_{\text{train}} \subseteq \mathcal{P}$,  tolerance $\varepsilon>0$, max. dimension $N_\text{max}\in\N$ 
	\State choose $\mu^{(1)}\in\cP_{\text{train}}$, compute snapshot $\zeta_1:=x_{\mu^{(1)}}^\cN$\newline and optimal test function $\eta_1$ with $\zeta_1=B^* \eta_1$
		\label{alg:WG:line1}
	\State \textbf{Initialize} $S_{1} \leftarrow \{ \mu^{(1)}\}$, 
	    $X_{1}:= \text{span} \{ \zeta_1\}$, $Y_{1}:= \text{span} \{ \eta_1\}$, $N:=1$
	\While{$N < N_\text{max}$}
		\State \textbf{if} $ \max\limits_{\mu \in \cP_{\text{train}}} \Delta_{N}(\mu) \leq \varepsilon$ \textbf{then} \Return
		\State $\mu^{(N+1)}\leftarrow \arg \max\limits_{\mu \in \cP_{\text{train}}}  \Delta_{N}(\mu)$
			\label{alg:mu max}
		\State compute snapshot $\zeta_{N+1}:=x_{\mu^{(N+1)}}^\cN$ and optimal test function $\eta_{N+1}$
		\State $S_{N+1} \leftarrow S_{N} \cup \{\mu^{(N+1)}\}$, 
		    $X_{N+1}:= X_N\oplus\text{span} \{ \zeta_{N+1}\}$,
		    $Y_{N+1}:= Y_N\oplus\text{span} \{\eta_{N+1}\}$
		\State $ N \leftarrow N+1$
	\EndWhile
    \Statex \textbf{output:} set of chosen parameters $S_{N}$, reduced spaces $X_N$, $Y_N$ 
    \end{algorithmic}
\end{algorithm}

\section{Numerical experiments}\label{Sec:4}
In this section, we report on results of some of our numerical experiments. Our main focus is on the numerical solution of the ultraweak form of the PDAE, the error estimation and the quantitative reduction. We solve the arising linear systems for the Petrov-Galerkin and the reduced system by \textsc{Matlab}'s backslash operator, see also our remarks in Section \ref{Sec:6} below. 
The codes for producing the subsequent results is available via  \url{https://github.com/mfeuerle/Ultraweak_PDAE}.

\subsection{Serial RLC circuit}
We start by a standard problem which (in some cases) admits a closed formula for the analytical solution. This allows us to monitor the exact error and a comparison with standard time-stepping methods. Our particular interest is the approximation property of the ultraweak approach, which is an $L_2$-approximation.

The serial RLC circuit consists of a resistor with resistance $R$, an inductor with inductance $L$, a capacitor with capacity $C$ and a voltage source fed by a voltage curve $f_{V_S}:\bar{I}\to\R$. Kirchhoff's circuit and further laws from electrical engineering yield a DAE with the data
\begin{gather*}
    E = \begin{pmatrix}
    1 & 0 & 0 & 0 \\
    0 & 1 & 0 & 0 \\
    0 & 0 & 0 & 0 \\
    0 & 0 & 0 & 0
    \end{pmatrix}, \;
    A = \begin{pmatrix}
    0 & 0 & L^{-1} & 0\\
	C^{-1} & 0 & 0 & 0\\
    R & 0 & 0 & -1\\
    0 & 1 & 1 & 1
    \end{pmatrix}, \;
    x = \begin{pmatrix}
    x_I \\ x_{V_C} \\ x_{V_L} \\ x_{V_R}
    \end{pmatrix}, \;
    f = \begin{pmatrix}
    0 \\ 0 \\ 0 \\ -f_{V_S}
    \end{pmatrix},
\end{gather*}
whose index is $k=1$. 
The solution $x$ consists of the electric current $x_I$ and the voltages at the capacitor $ x_{V_C}$, at the inductor $x_{V_L}$ and at the resistor $x_{V_R}$. 

\subsubsection*{Convergence of the Petrov-Galerkin scheme} 
\begin{figure}[!ht]
    \pgfplotstableread{ode_smooth.dat}{\odesmooth}
    \pgfplotstableread{uwf_smooth.dat}{\uwfsmooth}

    \begin{center}
        \begin{tikzpicture}[scale=0.53]
        \begin{axis}[
            xmin = 0, xmax = 12.6,
            xticklabel style={/pgf/number format/fixed},
            grid = both,
            minor tick num = 1,
            major grid style = {lightgray},
            minor grid style = {lightgray!25},
            width = \textwidth,
            height = 0.5\textwidth,
            legend cell align = {left},
            legend pos = north west,
            xlabel = {time $[s]$},
            ylabel = {voltage $V_L$ $[V]$}
        ]
            \addplot[black,line width=2.5pt, domain = 0:12.5,samples = 200,smooth] {-5*sin(deg(x))};
            \addplot[red,line width=2.5pt] table [x ={t}, y = {Vl}] {\uwfsmooth};
            \addplot[blue,dashed,line width=2.5pt] table [x = {t}, y = {Vl}] {\odesmooth};
            \legend{
                analytical solution,
                ultraweak,
                ode15i
            }
        \end{axis}
       \end{tikzpicture}
    \hfill
    \pgfplotstableread{ode_square.dat}{\odesquare}
    \pgfplotstableread{uwf_square.dat}{\uwfsquare}
    \begin{tikzpicture}[scale=0.53]
        \begin{axis}[
            xmin = 0, xmax = 0.105,
            xticklabel style={/pgf/number format/fixed},  
            scaled x ticks=true,
            grid = both,
            minor tick num = 1,
            major grid style = {lightgray},
            minor grid style = {lightgray!25},
            width = 0.9\textwidth,
            height = 0.5\textwidth,
            legend cell align = {left},
            legend pos = south east,
            xlabel = {time $[s]$},
            ylabel = {voltage $V_L$ $[V]$}
        ]
            \addplot[red,line width=2.5pt] table [x ={t}, y = {Vl}] {\uwfsquare};
            \addplot[blue,dashed,line width=2.5pt] table [x = {t}, y = {Vl}] {\odesquare};
            \legend{
                ultraweak,
                ode15i
            }
        \end{axis}
    \end{tikzpicture}   

    \caption{Serial RLC circuit, exact voltage at the inductor; comparison of time-stepping (ode15i -- blue) and ultraweak (red) approximation 
    for smooth $f_{V_S}^{\text{smooth}}$ (left, including analytical solution)  
    and discontinuous $f_{V_S}^{\text{disc}}$ (right) right-hand side.}
    \label{fig:rlc:ex_sol}
    \end{center}
\end{figure}
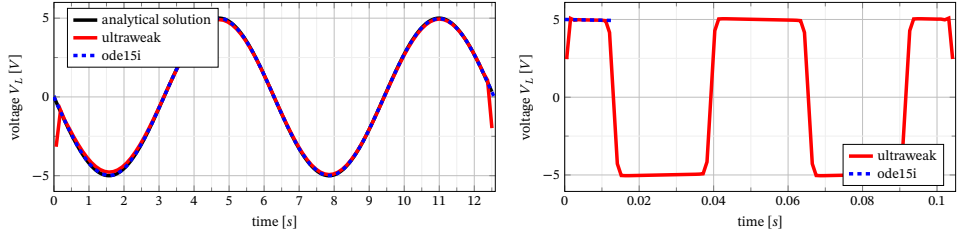

In Figure \ref{fig:rlc:ex_sol}, we compare the exact solution with approximations generated by a standard  time-stepping scheme (using \textsc{Matlab}'s fully implicit variable order solver with adaptive step size control \textit{ode15i}, \cite{MATLAB}) and by our ultraweak formulation from \S\ref{Sec:Disc}. We choose two specific examples for $f_{V_S}$, namely a smooth and a discontinuous one,
\begin{align*}
    f_{V_S}^{\text{smooth}}(t) := \sin\big(\frac{4\pi}{T}t\big),
    \qquad
    f_{V_S}^{\text{disc}}(t) := \operatorname{sign}
    \Big(\cos\big(\frac{4\pi}{T}t\big)\Big).
\end{align*}
For the smooth right-hand side (left graph in Figure \ref{fig:rlc:ex_sol}), both \textit{ode15i} and the ultraweak method give good results. Concerning the deviations for the ultraweak approach at the start and end time, we recall that the ultraweak form yields an approximation in $L_2$, so that pointwise comparisons are not necessarily meaningful.

In the discontinuous case, existence of a classical solution cannot be guaranteed by the above arguments. In particular, there is no closed solution formula. As we see in the right graph in Figure \ref{fig:rlc:ex_sol}, \textit{ode15i} stops at the first jump. This is to be expected, since $f_{V_S}^{\text{disc}}\not\in C^0(\bar{I})$, so that the solution lacks sufficient regularity to guarantee convergence of a time-stepping scheme like \textit{ode15i} (even though it is an adaptive variable order method). We could resolve the jumps even better by choosing more time steps $K$, while \textit{ode15i} still fails. We conclude that the ultraweak method also converges for problems lacking regularity.

\subsubsection*{Convergence rate}
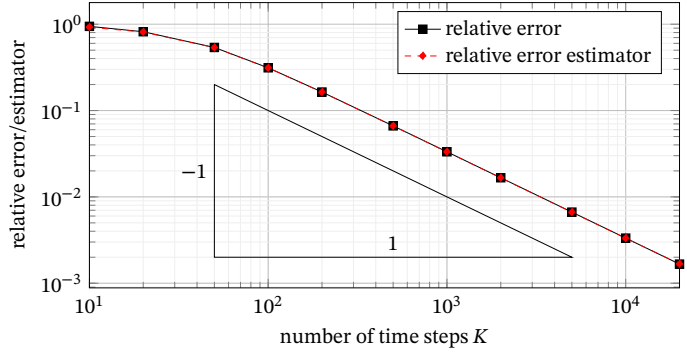
\begin{figure}[!htb]
    \pgfplotstableread{err_est_over_K.dat}{\errEstOverK}
    \begin{center}
    \begin{tikzpicture}[scale=0.8]
        \begin{axis}[
            xmode=log,
            ymode=log,
            xmin = 10, xmax = 20000,
            xticklabel style={/pgf/number format/fixed},
            grid = both,
            minor tick num = 1,
            major grid style = {lightgray},
            minor grid style = {lightgray!25},
            width = 0.9\textwidth,
            height = 0.5\textwidth,
            legend cell align = {left},
            legend pos = north east,
            xlabel = {number of time steps $K$},
            ylabel = {relative error/estimator}
        ]
            \addplot[black, mark = square*] table [x = {K}, y = {err}] {\errEstOverK};
            \addplot[red, mark = diamond*, dashed] table [x ={K}, y = {err_est}] {\errEstOverK};
            \draw (50,0.2) -- (5000,0.002) -- node [above] {$1$} (50,0.002) -- node [left] {$-1$} (50,0.2);
            \legend{
                relative error, 
                relative error estimator 
            }
        \end{axis}
    \end{tikzpicture}

   \caption{Relative error $\norm{x^*-x^\cN}_{L_2}/\norm{x^*}_{L_2}$ and relative error estimator $\Delta^\cN/\norm{x^*}_{L_2}$ w.r.t.\ to the analytical solution $x^*$ for increasing numbers of time steps $K$.}
    \label{fig:rlc:err_est_over_K}
    \end{center}
\end{figure}

Next, we investigate the rate of convergence for the ultraweak form. To that end, we use  $f_{V_S}^{\text{smooth}}$, since the analytical solution $x^*$ is known and we can thus compute the relative error $\norm{x^*-x^\cN}_{L_2}/\norm{x^*}_{L_2}$. Using the lowest order discretization as mentioned above (namely piecewise linear test functions $\psi_j$, which yield discontinuous trial functions $B^*\psi_i$), we can only hope for first order (w.r.t.\ the number of time steps $K$), which we see in Figure \ref{fig:rlc:err_est_over_K} and was observed in all cases we considered. We obtain higher order convergence by choosing test functions of higher order, provided the solution has sufficient smoothness.

Moreover, we compare the exact relative error with our error estimator (see \S\ref{Subsec:3.1}). Figure \ref{fig:rlc:err_est_over_K} shows a perfect matching confirming the error-residual identity \eqref{eq:res_identity} also for the numerically computed error estimator.

\subsection{Time dependent Stokes problem}
In order to investigate the quantitative performance of the model reduction, we consider a problem, which has often been used as a benchmark, \cite{MMESS,Schmidt20007,STYKEL2006,morwiki_stokes}, namely the time-dependent Stokes problem on the unit square $(0,1)^2$, discretized by a finite volume method on a uniform, staggered grid for the spatial variables with $n$ unknowns, \cite{morwiki_stokes}, where we choose $n=644$. The arising homogeneous fully linear DAE with output function $y : I\ra \R$ takes the form \eqref{eq:flDEA},
\begin{subequations}
\begin{align}
    E\dot{x}(t) - Ax(t) &= Bu(t)+g(t),\quad t\in I, \qquad x(0) = 0, \label{eq:stokes1}\\
    y(t) &= C x(t),
\end{align}
\end{subequations}
where $C\in\R^{1\times n}$ is an output matrix, $B\in\R^{n\times 1}$ is the control matrix, and $u:I\to\R$ is a control, which serves as a parameter $\mu \equiv u$ as described in \S\ref{Sec:flDAE} above.\footnote{We could also choose larger input/output-dimensions.} We use a parameter-independent initial condition, so that $g_\mu\equiv g$ and $Q_x=1$.

In order to combine system theoretic model reduction with the Reduced Basis Method from \S\ref{Sec:MOR}, we use the system theoretic model order reduction package \cite{MMESS}. In particular, we use  Balanced Truncation (BT) from \cite{STYKEL2006} during a preprocessing step to reduce the above system of dimension $n$ to a system
\begin{subequations}
\begin{align}
    \hat{E}\dot{\hat{x}}(t) - \hat{A}\hat{x}(t) 
    &= \hat{B}u(t) +\hat{g}(t),\quad t\in I,
    	\qquad \hat{x}(0) = 0,\label{eq:BTstokes1}\\
    y(t) &= \hat{C}\hat{x}(t),
\end{align}
\end{subequations}
with $\hat{E},\hat{A} \in \R^{\hat{n}\times\hat{n}}$, $\hat{B}\in \R^{\hat{n}\times 1}$, $\hat{C}\in \R^{1 \times \hat{n}}$ as well as $\hat{x}, \hat{g}: \bar{I} \ra \R^{\hat{n}}$ and $\hat{n}\ll n$. We note that the resulting reduced system typically provides regular matrices $\hat{E},\hat{A}$. Then, the reduced system is a linear time-invariant system (LTI), which is an easier problem than a DAE and in fact a special case. Hence, our presented approach is still valid, even though designed for PDAEs. For an ultraweak formulation of LTI systems, we refer to  \cite{Feuerle2019}.

\begin{remark}
We use the RBM here for deriving a certified reduced approximation of the state $x$. If we would want to control the output $y$ along with a corresponding error estimator $\Delta_N^y$, it is fairly standard in the theory of RBM to use a primal-dual approach with a second (dual) reduced basis, e.g. \cite{rozza2016,quarteroni2015}. For simplicity of exposition, we leave this to future work and compute the output from the state by $\hat{C}\hat{x}(t)$, resp.\ $Cx(t)$.\hfill$\diamond$
\end{remark}

\subsubsection*{Discretization of the control within the RBM} Since we use a variational approach, we are in principle free to choose any discretization for the control (we only need to compute inner products with the test basis functions). We tested piecewise linear discretizations as described in \S\ref{Sec:Disc} for different step sizes $K_u/T$, where $K_u$ might be different from $K$, which we choose for discretizing the state. 
Doing so, the parameter reads $\mu = \parens{u(t_0),\dots,u(t_{K_u})}^T \in \cP \equiv \R^{K_u+1}$, i.e., the parameter dimension is $P= K_u +1$, which might be large. Large parameter dimensions are potentially an issue for the RBM since the curse of dimension occurs. Hence, we investigate if we can reduce $K_u$ within the RBM.

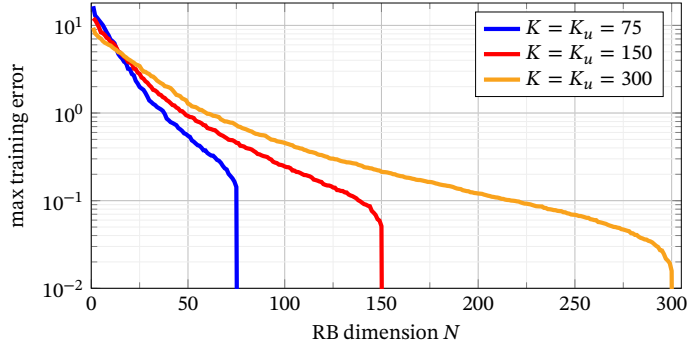
\begin{figure}[!ht]
    \begin{center}
    \pgfplotstableread{exponential_decay_greedy.dat}{\greedyDecay}

\begin{tikzpicture}[scale=0.8]
        \begin{axis}[
            ymode=log,
            xmin = 0, 
            xmax = 305,
            ymin = 0.01, 
            ymax = 18,
            xticklabel style={/pgf/number format/fixed},
            grid = both,
            minor tick num = 1,
            major grid style = {lightgray},
            minor grid style = {lightgray!25},
            width = 0.9\textwidth,
            height = 0.5\textwidth,
            legend cell align = {left},
            legend pos = north east,
            xlabel = {{RB dimension $N$}},
            ylabel = {max training error}
        ]
            \addplot[blue,line width=2pt] table [x = {N}, y = {greedy_err_K075}] {\greedyDecay};
            \addplot[red,line width=2pt] table [x = {N}, y = {greedy_err_K150}] {\greedyDecay};
            \addplot[YellowOrange,line width=2pt] table [x = {N}, y = {greedy_err_K300}] {\greedyDecay};
            \legend{
                $K=K_u=75$,
                $K=K_u=150$,
                $K=K_u=300$
            }
        \end{axis}
    \end{tikzpicture}

    \caption{Maximal greedy training error $\max_{\mu \in \cP_{\text{train}}}  \Delta_{N}(\mu)$ for different time resolutions $K_u=K\in\braces{75,150,300}$ over the reduced dimension $N$.}
    \label{fig:stokes:greedy_decay}
    \end{center}
\end{figure}

In order to answer this question, we apply Algorithm \ref{Alg:Greedy} to the time-dependent Stokes problem \eqref{eq:stokes1} (without BT) setting $\varepsilon = 0$, $N_\text{max} = Q_f$ from \eqref{eq:affinedecomp} (i.e., $Q_f = P+1$ for the fully linear system with parameter-independent initial value) and $P_\text{train}$ consisting of $500$ random vectors for $K_u\equiv K\in\braces{75,150,300}$, i.e., $\cN=48\,524,96\,824,193\,424$, where $d=224$. For these three cases, we investigate the \emph{max greedy training error}, i.e., $\max_{\mu \in \cP_{\text{train}}} \Delta_{N}(\mu)$. The results in Figure \crefl{fig:stokes:greedy_decay} show an exponential decay w.r.t.\ the dimension $N$ of the reduced system with slower decay as $K$ grows. This is to be expected as the discretized control space is much richer for growing $K_u$ and the reduced model has to be able to represent this variety. However, in relative terms (i.e., reduced size $N$ compared with full size $K$), we see that the compression rates are almost the same. This shows that the RBM can effectively reduce the system no matter how strong the influence of the control on the state is. It is expected that this potential is even more pronounced if a primal-dual RBM is used for the output. 

Next, we note that for $A \equiv A_\mu$ as \eqref{eq:stokes1}, the reduced model is always \emph{exact} for $N \geq Q_f$, which explains the drop off of the curves in Figure \crefl{fig:stokes:greedy_decay}. For fully linear DAEs, a reduced model with $N\geq Q_f = K_u+2$ is always exact. Hence, if $m\ll n$ (here $m=1 \ll 644=n$), we obtain an exact reduced model of dimension $N = Q_f = P+Q_x =m(K_u+1)+1 \ll nK +d =: \cN$. Even though this seems to be attractive for low-dimensional outputs, we stress the fact that the reduced dimension still depends on the temporal dimension $K_u$, which might be large. Hence, a combination of a possibly small discretization of the control and a RBM seems necessary.

Let us comment on the error decay of the RBM produced by the greedy method using the error estimator derived from the ultraweak formulation of the PDAE. We  obtain exponential decay of the error, which in fact shows the potential of the RBM. The question if a given PDAE permits a fast decay of the greedy RBM error is well-known to be linked to the decay of the Kolmogorov $N$-width, \cite{BennerSIAM2017,rozza2016,UrbanCIME2022}, which is a property only of the problem at hand. In other words, if a PDAE can be reduced w.r.t.\ time, the greedy method will detect this.

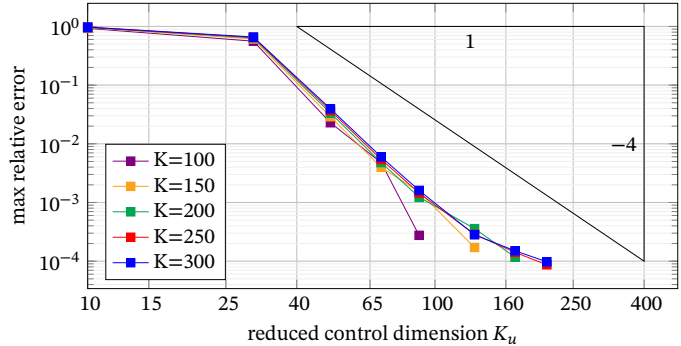
\begin{figure}[!ht]
    \begin{center}
    \pgfplotstableread{K100_n644_timereduction_sine_step.dat}{\timereductionKonehundred}
    \pgfplotstableread{K150_n644_timereduction_sine_step.dat}{\timereductionKonehundredfifty}
    \pgfplotstableread{K200_n644_timereduction_sine_step.dat}{\timereductionKtwohundred}
    \pgfplotstableread{K250_n644_timereduction_sine_step.dat}{\timereductionKtwohundredfifty}
    \pgfplotstableread{K300_n644_timereduction_sine_step.dat}{\timereductionKthreehundred}

    \pgfplotstableread{K100_n644_timereduction_sine_step.dat}{\timereductionKonehundred}
    \pgfplotstableread{K150_n644_timereduction_sine_step.dat}{\timereductionKonehundredfifty}
    \pgfplotstableread{K200_n644_timereduction_sine_step.dat}{\timereductionKtwohundred}
    \pgfplotstableread{K250_n644_timereduction_sine_step.dat}{\timereductionKtwohundredfifty}
    \pgfplotstableread{K300_n644_timereduction_sine_step.dat}{\timereductionKthreehundred}
    \hspace*{-5mm}
    \begin{tikzpicture}[scale=0.8]
        \begin{loglogaxis}[
            xmin = 10, 
            xmax = 500,
            xtick = {10, 15, 25, 40, 65, 100, 160, 250, 400},
            log x ticks with fixed point,
            grid = both,
            minor tick num = 1,
            major grid style = {lightgray},
            minor grid style = {lightgray!25},
            width = 0.9\textwidth,
            height = 0.5\textwidth,
            legend cell align = {left},
            legend pos = south west,
            xlabel = {{reduced control dimension $K_u$}},
            ylabel = {{max relative error}},
        ]
            \addplot[violet, mark = square*] table [x = {res}, y = {max_err}] {\timereductionKonehundred};
            \addplot[YellowOrange, mark = square*] table [x = {res}, y = {max_err}] {\timereductionKonehundredfifty};
            \addplot[Green, mark = square*] table [x ={res}, y = {max_err}] {\timereductionKtwohundred};
            \addplot[red, mark = square*] table [x ={res}, y = {max_err}] {\timereductionKtwohundredfifty};
            \addplot[blue, mark = square*] table [x ={res}, y = {max_err}] {\timereductionKthreehundred};
            \draw (40,1) -- (400,0.0001) -- node [left] {$-4$} (400,1) -- node [below] {$1$} (40,1);
            \legend{
                K=100, 
                K=150,
                K=200, 
                K=250,
                K=300,
            }
        \end{loglogaxis}
        \end{tikzpicture}      
        \caption{Max error for control dimensions of size $K_u < K$.}
        \label{fig:stokes:timereduction}
    \end{center}
\end{figure}

The results in Figure \ref{fig:stokes:greedy_decay} use $K_u=K$. The next question is how the error behaves for $K_u<K$. To this end, we determine the error in the state w.r.t.\ the full resolution, i.e., we compare the state derived from the control with $K_u$ degrees of freedom with the state of the fully resolved control. In Figure \ref{fig:stokes:timereduction}, we display errors for different values for $K$. We obtain fast convergence, which again shows the significant potential for a reduced temporal discretization of the control.

\subsubsection*{Combination with BT / RBM error decay}
Next, we wish to investigate if a combination of a system theoretic MOR (here BT) and an RBM-like reduction w.r.t.\ time can be combined. To this end, we fix the temporal resolution (i.e., the number of time steps, here $K=K_u=300$) and determine the RBM error using Algorithm \ref{Alg:Greedy} for the full and the BT-reduced system. We use \cite{MMESS} to compute the BT from \cite{STYKEL2006} and obtain a LTI system of dimension $\hat{n}=5$. 

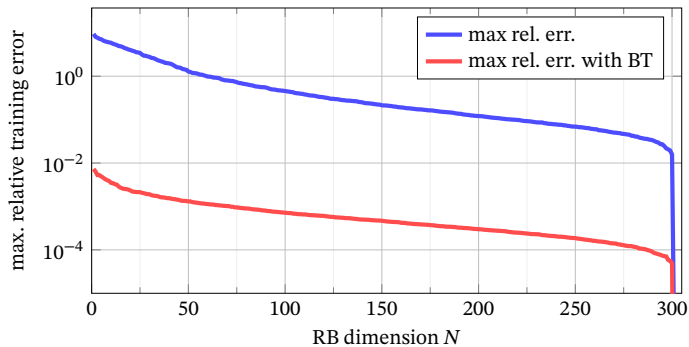
\begin{figure}[!ht]
    \begin{center}
    \pgfplotstableread{./greedy_compare_RB_BT.dat}{\RBBTcomparison}   

    \begin{tikzpicture}[scale=0.8]
        \begin{axis}[
            ymode=log,
            xmin = 0, 
            xmax = 305,
            ymin = 1e-5, 
            xticklabel style={/pgf/number format/fixed},
            grid = both,
            minor tick num = 1,
            major grid style = {lightgray},
            minor grid style = {lightgray!25},
            width = 0.9\textwidth,
            height = 0.5\textwidth,
            legend cell align = {left},
            legend pos = north east,
            xlabel = {RB dimension $N$},
            ylabel = {max.\ relative training error}
        ]
            \addplot[blue!70,line width=2pt] table [x ={N}, y = {RB}] {\RBBTcomparison};
            \addplot[red!70, line width=2pt] table [x ={N}, y = {RBBT}] {\RBBTcomparison};
            \legend{
                max rel.\ err., 
                max rel.\ err. with BT,
            }
        \end{axis}
    \end{tikzpicture}
    \caption{Maximal RBM relative error decay over the reduced dimension $N$ for the full system \crefl{eq:stokes1} (blue) and for the reduced system \crefl{eq:BTstokes1} with $K=300$ (red).
    }
    \label{fig:stokes:rb_bt_compare}
    \end{center}
\end{figure}

The results are shown in Figure \ref{fig:stokes:rb_bt_compare}, where we again show the maximal training error. As we see, the error for the BT-reduced system is smaller than the original one\footnote{This remains true even after additionally normalizing the training error by the dimensions of the DAE ($n$ and $\hat{n}$, respectively) or the dimension of the resulting linear system ($Kn+d$ and $K\hat{n}$, respectively).}, which in fact indicates that we can combine both methods. We get similar results for other choices of $K$. This shows that there is as much \enquote{reduction potential} in the reduced system \crefl{eq:BTstokes1} as in the original system \crefl{eq:stokes1}. In other words, a combination of BT and RBM shows significant compression potential.

\section{Conclusions and outlook}\label{Sec:6}
In this paper, we introduced a well-posed ultraweak formulation for DAEs and an  optimally stable Petrov-Galerkin discretization, which admits a sharp error bound. The scheme shows the expected order of convergence depending on the regularity of the solution and the smoothness of the trial functions. The scheme also converges in low-regularity cases, where classical standard time-stepping schemes fail. Moreover, the stability of the Petrov-Galerkin scheme allows us to choose \emph{any} temporal discretization without satisfying other stability criteria like a CFL condition.

Based upon the ultraweak framework, we introduced a model order reduction in terms of the Reduced Basis Method with an error/residual identity. We have obtained fast convergence and the possibility to combine the RBM for a reduction w.r.t.\ time with system theoretic methods such as Balanced Truncation to reduce the size of the system.

There are several open issues for future research. We already mentioned a primal-dual RBM for an efficient reduction of the output, the generalization to parameter-dependent matrices $A_\mu$ and more general DAEs (not only fully linear). We also mentioned that the system matrix is a sum of Kronecker products of high dimension, which calls for specific solvers as in \cite{Henning2021} for the (parameterized) wave equation. Another issue in that direction is the need for a basis of $\kernel(E^T)$, which might be an issue for high-dimensional problems.

\printbibliography
\end{document}